\documentclass[11pt]{article}
\usepackage[letterpaper,margin=1in]{geometry}
\usepackage[T1]{fontenc}
\usepackage{newtxtext}

\usepackage{amsthm,amsmath,amsfonts,amssymb}
\usepackage[numbers,sort&compress]{natbib}
\usepackage{graphicx}
\usepackage{booktabs}
\usepackage[hidelinks]{hyperref}

\theoremstyle{plain}
\newtheorem{theorem}{Theorem}

\newtheorem{proposition}{Proposition}
\newtheorem{corollary}{Corollary}
\theoremstyle{definition}
\newtheorem{assumption}{Assumption}
\DeclareMathOperator{\E}{E}
\DeclareMathOperator{\var}{var}

\DeclareMathOperator{\pr}{pr}
\DeclareMathOperator{\sgn}{sgn}
\newcommand{\dd}{\,\mathrm{d}}
\newcommand{\dto}{\overset{d}{\longrightarrow}}
\newcommand{\pto}{\overset{p}{\longrightarrow}}

\newcommand{\RS}{\rho}
\newcommand{\TK}{\tau}
\newcommand{\rhon}{\widehat{\rho}}
\newcommand{\taun}{\widehat{\tau}}

\hypersetup{pdftitle={Nonsense rank correlations between independent random walks},pdfauthor={Zhongming Xie and Liping Zhu},pdfsubject={Preprint}}

\title{Nonsense rank correlations between independent random walks}
\author{Zhongming Xie \qquad Liping Zhu\\[0.5ex]
\small Institute of Statistics and Big Data\\
\small Renmin University of China\\
\small \texttt{xzm1@ruc.edu.cn} \qquad \texttt{zhu.liping@ruc.edu.cn}}
\date{}

\begin{document}
\maketitle
\begin{abstract}
Yule's nonsense correlation persists for Spearman's rho and Kendall's tau between independent random walks. Neither coefficient concentrates at zero as the sample size grows. For centred random walks with finite increment variances, we derive a universal joint limit expressed through Brownian occupation functionals. Integral formulas give limiting variances of approximately 0.222 and 0.108, respectively, below the Pearson value of 0.241. Rank correlations thus reduce the limiting dispersion on their respective coefficient scales. We prove that both marginal limits have symmetric densities and represent them exactly as limits of polynomial densities determined by Gaussian moment integrals. Spearman's conventional $t$ test and Kendall's normal approximation are asymptotically invalid, with rejection probabilities tending to one under process independence.
\end{abstract}

\section{Introduction}

Correlations between independent time series can remain widely dispersed as the sample size increases. The experiments of Yule~\cite{yule1926} showed that independent random walks can produce Pearson correlations much larger than independent-observation theory suggests. Ernst et al.~\cite{ernst2017} derived an integral representation for the second moment of the continuous-time Pearson correlation between independent Brownian motions. Its approximate value, 0.241, quantifies the dispersion that remains under process independence. Subsequent work studied finite Gaussian random walks \citep{ernst2023} and the moments and density of the Brownian correlation \citep{ernst2025}.

Sz\'ekely and Rizzo~\cite[\S\,21.1]{szekelyrizzo2023} attribute nonsense correlation to serial dependence and nonstationarity. They describe the phenomenon as ``contagious'' and illustrate it for distance correlation between independent random walks. Zhu and Xie~\cite{zhuxie2026} establish a nondegenerate Brownian occupation limit for sample distance correlation and show that permutation calibration rejects process independence with probability tending to one.

Random-walk levels are serially dependent, and their marginal variances increase with time. The resulting variation in marginal scale motivates using Spearman's rho and Kendall's tau, which depend only on relative order. We ask whether discarding observation magnitudes can attenuate nonsense correlation between independent random walks.

In related work on rank-based independence measures, Zhou et al.~\cite{zhou2024} show that Hoeffding's $D$, Blum--Kiefer--Rosenblatt's $R$ and Bergsma--Dassios--Yanagimoto's $\tau^*$ are strictly increasing functions of the absolute Pearson correlation for bivariate normal variables. For time series, Schaworonkow et al.~\cite{schaworonkow2015} simulated pairs of mutually independent series, each with power-law temporal dependence. They found that Spearman's sample correlation could be more dispersed than predicted by independent-observation theory. For independent centred random walks, we give a limiting-distribution explanation of this phenomenon and show that it also occurs for Kendall's tau.

Our results have three parts. First, we establish the joint limiting distribution of Spearman's sample correlation $\widehat{\rho}$ and Kendall's sample correlation $\widehat{\tau}$ under the stated increment assumptions. We derive integral formulas for their limiting variances, which are approximately $0.222$ and $0.108$, respectively. Both are below the Pearson value of $0.241$ obtained by Ernst et al.~\cite{ernst2017}. These smaller variances quantify attenuation of nonsense correlation on the respective coefficient scales, although neither coefficient concentrates at zero. Second, we obtain exact $L^1$-convergent representations of the two marginal limiting densities. Third, we show that Spearman's conventional $t$ test and Kendall's normal approximation reject independence of the two random walks with probability tending to one at every fixed nominal significance level.

Section~\ref{sec:occupation} establishes the joint Brownian limit. Section~\ref{sec:moments} derives the limiting variances and endpoint probabilities, and Section~\ref{sec:density} gives the marginal densities. Section~\ref{sec:calibration} explains the failure of conventional rank tests and discusses Brownian calibration and power. Section~\ref{sec:numerical} reports numerical integration and simulation results, and Section~\ref{sec:discussion} discusses implications and further questions. Proofs and computational details appear in the Supplement~\cite{xiezhu2026supp}.

\section{The occupation-measure limit}\label{sec:occupation}
\subsection{Model and statistics}

Consider the levels of two random walks,
\begin{equation}\label{eq:model}
X_k=\sum_{i=1}^k\varepsilon_i,\qquad Y_k=\sum_{i=1}^k\eta_i,\qquad X_0=Y_0=0.
\end{equation}
\begin{assumption}[Random walks]\label{ass:random-walk}
The increment sequences in \eqref{eq:model} satisfy the following conditions.
\renewcommand{\labelenumi}{(\roman{enumi})}
\begin{enumerate}
\item The sequences $(\varepsilon_i)_{i\geq1}$ and $(\eta_i)_{i\geq1}$ are independent of each other.
\item Within each sequence, the increments are independent and identically distributed, with a continuous distribution.
\item The increments satisfy $\E(\varepsilon_1)=\E(\eta_1)=0$.
\item The variances satisfy $\var(\varepsilon_1)=\sigma_X^2\in(0,\infty)$ and $\var(\eta_1)=\sigma_Y^2\in(0,\infty)$.
\end{enumerate}
\end{assumption}

The two increment distributions need not be equal, symmetric or Gaussian. Continuity ensures that each sample of path levels has no ties, almost surely. The variables $X_k$ and $Y_l$ are independent for every $k,l$. Successive values within each walk are dependent because each value includes all the increments accumulated at earlier times.

For $n\geq2$ and $k=1,\ldots,n$, let $R_k$ and $S_k$ be the ranks of $X_k$ and $Y_k$ among their respective samples of size $n$. We study Spearman's rho and Kendall's tau, defined by
\begin{align}
\rhon&=\frac{12}{n(n^2-1)}\sum_{k=1}^n\{R_k-(n+1)/2\}\{S_k-(n+1)/2\},\label{eq:spearman-sample}\\
\taun&=\frac{2}{n(n-1)}\sum_{1\leq i<j\leq n}\sgn(X_i-X_j)\sgn(Y_i-Y_j),\label{eq:kendall-sample}
\end{align}
where $\sgn(0)=0$. Both statistics are unchanged by separate strictly increasing transformations of the two coordinates. In particular, dividing the walks by $\sigma_X n^{1/2}$ and $\sigma_Y n^{1/2}$ leaves their ranks and pairwise signs unchanged. This allows us to use the usual Brownian scaling of the paths without rescaling the rank statistics themselves.

\subsection{Occupation-measure functionals}

For continuous paths $f,g$ on $[0,1]$, let $\mu_{f,g}(\mathcal A)$ denote the proportion of time for which the pair $(f(t),g(t))$ lies in $\mathcal A$, given by
$$
\mu_{f,g}(\mathcal A)=\int_0^1 I\{(f(t),g(t))\in\mathcal A\}\dd t
$$
for Borel sets $\mathcal A\subseteq\mathbb R^2$. Equivalently, $\mu_{f,g}$ is the distribution of $(f(V),g(V))$ when $V$ is uniform on $[0,1]$ and the paths are held fixed.

Write $F_f(x)=\int_0^1 I\{f(t)\leq x\}\dd t$ for the first marginal distribution function, and define $F_g$ similarly. The quantity $F_f\{f(t)\}$ is the proportion of times at which the path is at or below its level at time $t$. It is the continuous-time counterpart of the normalized sample rank $R_k/n$. If $F_f$ is continuous, $F_f\{f(V)\}$ is uniform on $[0,1]$ under the random choice of $V$. Consequently,
\begin{equation}\label{eq:uniform-occupation}
\int_0^1 F_f\{f(t)\}\dd t=\frac12,\qquad \int_0^1 [F_f\{f(t)\}-1/2]^2\dd t=\frac1{12}.
\end{equation}
The same identities hold for $g$ when $F_g$ is continuous. They give the fixed centring and scaling constants in the rank functionals defined by
\begin{align}
\mathcal R(f,g)&=12\int_0^1 [F_f\{f(t)\}-1/2][F_g\{g(t)\}-1/2]\dd t,\label{eq:R-functional}\\
\mathcal T(f,g)&=\int_0^1\!\int_0^1\sgn\{f(s)-f(t)\}\sgn\{g(s)-g(t)\}\dd s\dd t.\label{eq:T-functional}
\end{align}
When the marginal distribution functions $F_f$ and $F_g$ are continuous, these are the population Spearman and Kendall correlations of $(f(V),g(V))$. Both lie in $[-1,1]$. The first averages products of centred normalized ranks at a common time. The second averages concordance signs over two independently and uniformly chosen times.

For random paths, the occupation measure describes one realized pair of trajectories sampled at a common uniform time. Independently generated paths can have similar or opposite ordering patterns over time.

\subsection{A universal joint limit}

The normalized ranks $R_k/n$ and $S_k/n$ are the proportions of observed levels at or below $X_k$ and $Y_k$ within their respective series. Their continuous-time counterparts are $F_f\{f(t)\}$ and $F_g\{g(t)\}$. Spearman's sum over observation times corresponds to the integral in \eqref{eq:R-functional}, while Kendall's sum over pairs of times corresponds to the double integral in \eqref{eq:T-functional}. The following theorem identifies the limits of the sample coefficients by evaluating these functionals at independent Brownian paths.

\begin{theorem}\label{thm:limit}
Under Assumption~\ref{ass:random-walk}, let $W_1,W_2$ be independent standard Brownian motions on $[0,1]$. Then, as $n\longrightarrow\infty$,
\begin{equation}\label{eq:joint-limit}
(\rhon,\taun)\dto(\RS,\TK),\qquad \RS=\mathcal R(W_1,W_2),\quad \TK=\mathcal T(W_1,W_2).
\end{equation}
For all nonnegative integers $a,b$, the mixed moments also satisfy
$$
\E(\rhon^a\taun^b)\longrightarrow\E(\RS^a\TK^b).
$$
\end{theorem}

Under Assumption~\ref{ass:random-walk}, the limiting joint law is the same for all choices of the two increment distributions and their variances. The two sample correlations vary across independently generated pairs of walks of length $n$. As $n$ increases, their joint distribution approaches that of $(\RS,\TK)$. Their means, variances and covariance also converge to those of this limiting pair.

\section{Variances and endpoint behaviour}\label{sec:moments}
\subsection{A Gaussian increment kernel}

Reflection of one Brownian path changes the sign of either coefficient without changing its distribution. Both means are therefore zero, so the second moments computed below equal the limiting variances.

Define, for $s\neq t$ and $u\neq v$,
\begin{equation}\label{eq:q}
q(s,t;u,v)=\frac{\min(s,u)-\min(s,v)-\min(t,u)+\min(t,v)}{(|s-t|\,|u-v|)^{1/2}}.
\end{equation}
This is the correlation between $W(s)-W(t)$ and $W(u)-W(v)$ within one standard Brownian motion. We set $q$ to zero on $s=t$ or $u=v$, which does not affect the integrals below.

The following Gaussian sign identity is used to derive both limiting variances. It gives the covariance kernel in \eqref{eq:C} for Spearman's rho and the squared-arcsine integrand in \eqref{eq:tau-variance} for Kendall's tau. If $G,H$ are centred jointly Gaussian variables with unit variances and correlation $r$, then
\begin{equation}\label{eq:arcsine}
\E\{\sgn(G)\sgn(H)\}=\frac2\pi\arcsin(r).
\end{equation}
We apply the identity to two standardized increments from the same Brownian motion, using $q$ as their correlation. A short proof of the identity is included in Section S3 of the Supplement~\cite{xiezhu2026supp}.

\subsection{Spearman's limiting variance}

Set $A_W(t)=F_W\{W(t)\}-1/2$, the centred normalized rank of a Brownian path at time $t$. Brownian reflection symmetry gives $\E\{A_W(t)\}=0$ for each fixed $t$. Thus $C(s,t)=\E\{A_W(s)A_W(t)\}$ is the covariance kernel of the rank process along one path.

\begin{theorem}\label{thm:spearman-moment}
The covariance kernel and limiting variance satisfy
\begin{align}
C(s,t)&=\frac1{2\pi}\int_0^1\!\int_0^1\arcsin\{q(s,u;t,v)\}\dd u\dd v,\label{eq:C}\\
\var(\RS)&=144\int_0^1\!\int_0^1 C(s,t)^2\dd s\dd t\approx0.222.\label{eq:rho-variance}
\end{align}
\end{theorem}

The variance in Theorem~\ref{thm:spearman-moment} corresponds to a limiting standard deviation of approximately $0.471$. Section~\ref{sec:quadrature} describes the numerical evaluation.

\subsection{Kendall's limiting variance}

We next quantify the limiting dispersion of Kendall's tau between independent random walks. The following theorem expresses its variance through the same Brownian increment correlation $q$ used for Spearman's rho.

\begin{theorem}\label{thm:kendall-moment}
The limiting Kendall variance satisfies
\begin{equation}\label{eq:tau-variance}
\var(\TK)=\frac4{\pi^2}\int_{[0,1]^4}\arcsin^2\{q(s,t;u,v)\}\dd s\dd t\dd u\dd v\approx0.108.
\end{equation}
\end{theorem}

The variance in Theorem~\ref{thm:kendall-moment} corresponds to a limiting standard deviation of approximately $0.328$. The two-dimensional representation in equation (S2) of the Supplement~\cite{xiezhu2026supp} is used for the numerical evaluation in Section~\ref{sec:quadrature}.

For comparison, Ernst et al.~\cite{ernst2017} obtain a variance of approximately $0.241$ for the Brownian Pearson correlation. Spearman's limiting variance of $0.222$ is slightly smaller, whereas Kendall's variance of $0.108$ is less than half as large.

\subsection{Reflection symmetry and endpoint probabilities}

The positive variances in Theorems~\ref{thm:spearman-moment} and~\ref{thm:kendall-moment}, together with Theorem~\ref{thm:limit}, rule out convergence in probability to zero. The next theorem gives a stronger joint property, showing that both coefficients can simultaneously be close to one or simultaneously close to minus one.

\begin{theorem}\label{thm:nondegenerate}
The limiting pair has the reflection symmetry $(\RS,\TK)\overset{d}{=}(-\RS,-\TK)$. In particular, both expectations are zero and every mixed moment of odd total degree vanishes. For $0<\delta<1$,
\begin{align}
\pr(\RS>\delta,\TK>\delta)&>0,\label{eq:positive-endpoint}\\
\pr(\RS<-\delta,\TK<-\delta)&>0.\label{eq:negative-endpoint}
\end{align}
\end{theorem}

Brownian scaling gives a complementary interpretation of the persistent dispersion. Rescaling time from $[0,L]$ to $[0,1]$ and dividing each path by $L^{1/2}$ gives the same Brownian law. Positive scaling preserves ranks, so the rank functionals defined with normalized time measure $\dd t/L$ have the same joint distribution for every $L>0$. A longer Brownian observation horizon therefore leaves their dispersion unchanged. Theorem~\ref{thm:limit} gives the corresponding large-sample conclusion for discrete random walks.

\section{Marginal densities of the Brownian rank limits}\label{sec:density}

We first show that the Brownian rank limits have marginal densities. We then represent these densities using Gaussian moment integrals.

\begin{proposition}\label{prop:density-existence}
The limiting variables $\RS$ and $\TK$ have marginal densities on $(-1,1)$ that are symmetric about zero, with no probability mass at the endpoints. In particular, $\pr(\RS=0)=\pr(\TK=0)=0$.
\end{proposition}

The absence of mass at zero means that the probability of either limit falling in a shrinking interval around zero tends to zero. Combined with Theorem~\ref{thm:limit}, this explains why the conventional tests in Section~\ref{sec:calibration}, whose correlation thresholds shrink to zero, reject independent random walks with probability tending to one.

The proof is given in Section S6 of the Supplement~\cite{xiezhu2026supp}.

\subsection{Exact moments as Gaussian integrals}

Both rank functionals can be expressed through products of signs of Brownian increments. Taking powers and integrating therefore reduces their moments to Gaussian sign expectations, extending the second-moment calculations in Section~\ref{sec:moments}. For two $k$-dimensional time vectors $\boldsymbol{s},\boldsymbol{t}$, let $\boldsymbol{\Sigma}(\boldsymbol{s},\boldsymbol{t})$ have entries $q(s_i,t_i;s_j,t_j)$, with $q$ defined in \eqref{eq:q}. Define
\begin{equation}\label{eq:gaussian-sign-all}
H_k(\boldsymbol{s},\boldsymbol{t})=(2\pi)^{-k/2}|\boldsymbol{\Sigma}(\boldsymbol{s},\boldsymbol{t})|^{-1/2}\int_{\mathbb R^k}\left\{\prod_{j=1}^k\sgn(z_j)\right\}\exp\left\{-(1/2)\boldsymbol z^{\mbox{\normalfont\tiny T}}\boldsymbol{\Sigma}(\boldsymbol{s},\boldsymbol{t})^{-1}\boldsymbol z\right\}\dd\boldsymbol z.
\end{equation}
The matrix $\boldsymbol{\Sigma}(\boldsymbol{s},\boldsymbol{t})$ is positive definite for almost every time configuration used below. At singular configurations, $H_k$ denotes the corresponding Gaussian sign expectation, which is also the continuous Gaussian-distribution limit when no individual increment has zero variance. Values on zero-length increments can be assigned arbitrarily without changing the time integrals. Set $H_0=1$.

For either $Q\in\{\RS,\TK\}$, the zeroth moment is one. Independence of the Brownian paths gives, for every integer $j\geq1$,
\begin{align}
\E(\TK^j)&=\int_{[0,1]^{2j}}H_j(\boldsymbol{s},\boldsymbol{t})^2\dd\boldsymbol{s}\dd\boldsymbol{t},\label{eq:tau-all-moments}\\
\E(\RS^j)&=3^j\int_{[0,1]^{3j}}H_j(\boldsymbol{s},\boldsymbol{u})H_j(\boldsymbol{s},\boldsymbol{v})\dd\boldsymbol{s}\dd\boldsymbol{u}\dd\boldsymbol{v}.\label{eq:rho-all-moments}
\end{align}
Odd-order Gaussian sign moments vanish, in agreement with Theorem~\ref{thm:nondegenerate}. At order two, these identities recover the variance formulas in Theorems~\ref{thm:spearman-moment} and~\ref{thm:kendall-moment}. At higher orders, they determine the coefficients needed for the density representation. Each integral is finite-dimensional, although its dimension increases with the moment order.

\subsection{A convergent density representation}

We use the moments to construct polynomial densities that converge to the marginal density in $L^1$. Let $b_{m,k}(u)=\binom{m}{k}u^k(1-u)^{m-k}$ for $0\leq u\leq1$ and $k=0,\ldots,m$. Define the polynomial coefficients $d_j(m,k)$ by
$$
(1+z)^k(1-z)^{m-k}=\sum_{j=0}^m d_j(m,k)z^j,
$$
and set
\begin{equation}\label{eq:density-weights}
w_{m,k}(Q)=2^{-m}\binom{m}{k}\sum_{j=0}^m d_j(m,k)\E(Q^j).
\end{equation}
The weights are expectations of the Bernstein basis polynomials evaluated at $(Q+1)/2$, as shown in Section S7 of the Supplement~\cite{xiezhu2026supp}. Equation~\eqref{eq:density-weights} determines them from the Gaussian moment integrals.

\begin{theorem}\label{thm:density-representation}
For $Q=\RS$ or $Q=\TK$, its marginal density has the exact representation
\begin{equation}\label{eq:exact-density}
f_Q=\lim_{m\longrightarrow\infty}p_{m,Q}\quad\hbox{in }L^1(-1,1),\qquad
p_{m,Q}(x)=\frac{m+1}{2}\sum_{k=0}^m w_{m,k}(Q)b_{m,k}\left(\frac{x+1}{2}\right).
\end{equation}
For every $m$, the weights are nonnegative and sum to one, and $p_{m,Q}$ is a symmetric probability density on $[-1,1]$.
\end{theorem}

For $m=1,2,3$, the polynomial densities depend only on the second moments, since symmetry makes all odd moments zero. For $Q\in\{\RS,\TK\}$ and $-1\leq x\leq1$, expanding \eqref{eq:exact-density} gives
\begin{align*}
p_{1,Q}(x)&=\frac12,\\
p_{2,Q}(x)&=\frac12+\frac{3\E(Q^2)-1}{16}(3x^2-1),\\
p_{3,Q}(x)&=\frac12+\frac{3\E(Q^2)-1}{8}(3x^2-1).
\end{align*}
The fourth moment is also needed for $m=4,5$, and the sixth moment enters when $m=6$. Section S11 of the Supplement~\cite{xiezhu2026supp} gives the exact expansions for these orders. Substituting the numerical moments from Section~\ref{sec:quadrature} gives the coefficients in Table~\ref{tab:polynomial-densities}.

\begin{table}[htbp]
\caption{Coefficients of the indicated powers of $x$ in $p_{m,Q}(x)$ on $[-1,1]$. Entries use the numerical moments and are rounded to three decimal places.}\label{tab:polynomial-densities}
\centering
\begin{tabular}{cccccc}
\toprule
$m$ & $Q$ & $1$ & $x^2$ & $x^4$ & $x^6$\\
\midrule
1 & $\RS,\TK$ & $0.500$ & $0.000$ & $0.000$ & $0.000$\\
2 & $\RS$ & $0.521$ & $-0.063$ & $0.000$ & $0.000$\\
2 & $\TK$ & $0.542$ & $-0.127$ & $0.000$ & $0.000$\\
3 & $\RS$ & $0.542$ & $-0.126$ & $0.000$ & $0.000$\\
3 & $\TK$ & $0.585$ & $-0.254$ & $0.000$ & $0.000$\\
4 & $\RS$ & $0.559$ & $-0.174$ & $-0.006$ & $0.000$\\
4 & $\TK$ & $0.622$ & $-0.374$ & $0.013$ & $0.000$\\
5 & $\RS$ & $0.573$ & $-0.208$ & $-0.019$ & $0.000$\\
5 & $\TK$ & $0.654$ & $-0.487$ & $0.038$ & $0.000$\\
6 & $\RS$ & $0.584$ & $-0.231$ & $-0.036$ & $0.000$\\
6 & $\TK$ & $0.683$ & $-0.592$ & $0.073$ & $0.001$\\
\bottomrule
\end{tabular}
\end{table}

\section{Conventional rank tests and their failure}\label{sec:calibration}

Spearman's $t$ test and Kendall's normal approximation use reference distributions for independent and identically distributed observation pairs with independent continuous coordinates. Their correlation thresholds tend to zero, while the coefficients from independent random walks converge to Brownian limits with positive variances. Throughout this section, the tests are two-sided with fixed nominal level $\alpha\in(0,1)$.

\subsection{Spearman's \texorpdfstring{$t$}{t} test}\label{sec:spearman-test}

For $n\geq3$, the usual approximate Spearman test uses the transformed statistic
$$
T_\rho=\rhon\left\{\frac{n-2}{1-\rhon^2}\right\}^{1/2}
$$
and a Student $t$ reference distribution with $n-2$ degrees of freedom \citep{hollander1973}. Write $t_{\nu,p}$ for the $p$th quantile of a Student $t$ distribution with $\nu$ degrees of freedom. The test rejects when $|\rhon|=1$ or $|T_\rho|>t_{n-2,1-\alpha/2}$. On the correlation scale, this rule is $|\rhon|>c_\rho$, where
\begin{equation}\label{eq:spearman-critical}
c_\rho=\frac{t_{n-2,1-\alpha/2}}{\{n-2+t_{n-2,1-\alpha/2}^2\}^{1/2}}=O(n^{-1/2}).
\end{equation}
\subsection{Kendall's normal approximation}\label{sec:kendall-test}

Under the same independent-observation null model, Kendall's coefficient has mean zero and exact variance
$$
\var(\taun)=\frac{2(2n+5)}{9n(n-1)}.
$$
The usual large-sample test compares the standardized statistic
$$
Z_\tau=\left\{\frac{9n(n-1)}{2(2n+5)}\right\}^{1/2}\taun
$$
with a standard normal distribution. Writing $z_p$ for its $p$th quantile, the test rejects when $|Z_\tau|>z_{1-\alpha/2}$. Equivalently, it rejects when $|\taun|>c_\tau$, where
\begin{equation}\label{eq:kendall-critical}
c_\tau=z_{1-\alpha/2}\left\{\frac{2(2n+5)}{9n(n-1)}\right\}^{1/2}=O(n^{-1/2}).
\end{equation}
Both conventional tests therefore use correlation thresholds of order $n^{-1/2}$.

\subsection{Why both tests fail for random-walk levels}\label{sec:test-failure}

Under Assumption~\ref{ass:random-walk}, the two walks are independent of one another, but successive levels within each walk are dependent. Theorem~\ref{thm:limit} gives $\rhon\dto\RS$ and $\taun\dto\TK$ without rescaling. The limiting standard deviations are approximately $0.471$ and $0.328$, respectively, by Theorems~\ref{thm:spearman-moment} and~\ref{thm:kendall-moment}. The sample correlations therefore remain dispersed as $n$ increases.

At $n=2000$ and nominal level $5\%$, the two conventional tests reject when the absolute correlations exceed $0.044$ and $0.029$, respectively. As $n$ grows, these thresholds approach zero while the null correlations continue to vary on the scales described above.

To state the consequence for both coefficients, let $(\widehat{Q},Q)$ denote either $(\rhon,\RS)$ or $(\taun,\TK)$. The following result applies to any nonnegative threshold tending to zero, including thresholds estimated from the data.

\begin{corollary}\label{cor:shrinking}
Under Assumption~\ref{ass:random-walk}, if a nonnegative critical value $c_n$ satisfies $c_n\pto0$, then
\begin{equation}\label{eq:shrinking-critical}
\pr(|\widehat{Q}|>c_n)\longrightarrow1.
\end{equation}
\end{corollary}

Both tests reject whenever the coefficient falls outside a shrinking interval around zero. Proposition~\ref{prop:density-existence} gives $\pr(Q=0)=0$, so the probability of remaining inside this interval tends to zero. Applying Corollary~\ref{cor:shrinking} to \eqref{eq:spearman-critical} and \eqref{eq:kendall-critical} therefore shows that both tests reject independent random walks with probability tending to one at every fixed nominal level. The proof is given in Section S8 of the Supplement~\cite{xiezhu2026supp}.

\subsection{Other tests based on independent observations}\label{sec:other-tests}

The same argument covers Spearman's normal approximation, which compares $Z_\rho=(n-1)^{1/2}\rhon$ with a standard normal distribution. Its correlation threshold is $z_{1-\alpha/2}/(n-1)^{1/2}$. A Fisher-type rule $|\operatorname{arctanh}(\widehat{Q})|>b_n$ is also covered whenever $b_n\pto0$, because its threshold on the correlation scale is $\tanh(b_n)$.

Exact rank tests and unrestricted permutation tests use the independent-observation rank distribution. Its variance is $1/(n-1)$ for Spearman and $2(2n+5)/\{9n(n-1)\}$ for Kendall, so their critical values are $O(n^{-1/2})$ (Section S8 of the Supplement~\cite{xiezhu2026supp}). Thus Corollary~\ref{cor:shrinking} applies to these tests as well. Unrestricted permutation of one walk's observations disrupts its temporal dependence, consistent with the distance-correlation result of Zhu and Xie~\cite{zhuxie2026}.

Resampling individual observations also yields shrinking thresholds. The moment bounds in Lemma S3 of the Supplement~\cite{xiezhu2026supp} cover Wald or centred-quantile tests from the paired bootstrap, tests using reference quantiles from independently resampled marginals, and Wald tests based on the ordinary delete-one jackknife variance. The independent-multiplier references with uniformly bounded scores described in Section S9 of the Supplement have variance $O(n^{-1})$ and hence shrinking critical values.

\subsection{Brownian calibration and power}\label{sec:brownian-calibration}

The Brownian limit supplies a reference distribution appropriate to the random-walk model. Choose a fixed $c\in(0,1)$ satisfying $\pr(|Q|>c)=\alpha$, and reject when $|\widehat{Q}|>c$. Such a threshold exists because Proposition~\ref{prop:density-existence} ensures that the tail probability varies continuously from one at $c=0$ to zero at $c=1$. This proposition also gives $\pr(Q=c)=\pr(Q=-c)=0$, so the rejection region is a continuity set of the limiting law. Under Assumption~\ref{ass:random-walk}, Theorem~\ref{thm:limit} therefore yields
\begin{equation}\label{eq:calibration}
\pr(|\widehat{Q}|>c)\longrightarrow\pr(|Q|>c)=\alpha.
\end{equation}
The threshold remains fixed as $n$ grows, matching the persistent dispersion of the rank coefficients. Section~\ref{sec:numerical} evaluates simulated approximations to these Brownian thresholds.

Correct size does not ensure consistency against dependent alternatives. Consider the fixed alternative in which the walks in \eqref{eq:model} have increments
\begin{equation}\label{eq:power-alternative}
\varepsilon_i=Z_i,\qquad \eta_i=(Z_i^2-1)/2^{1/2},
\end{equation}
where the $Z_i$ are independent standard normal variables. Each $\eta_i$ is a nonconstant function of $\varepsilon_i$, so the walks are dependent. Their increment covariance is zero, and their rank correlations have the same joint limit as under independence. Consequently, both tests using the fixed null critical values defined above have power converging to $\alpha$ at this alternative. Section S10 of the Supplement~\cite{xiezhu2026supp} proves this claim.

\section{Numerical results}\label{sec:numerical}

Numerical integration evaluates the limiting moment formulas and supplies the coefficients of the polynomial densities in Table~\ref{tab:polynomial-densities}. Path simulations examine whether the correlations remain dispersed at finite sample sizes and compare null rejection frequencies under conventional and approximate Brownian calibration.

\subsection{Numerical evaluation of the limiting moments}\label{sec:quadrature}

We use product Gauss--Legendre quadrature for Spearman's variance in \eqref{eq:rho-variance} and adaptive quadrature for Kendall's variance using equation (S2) of the Supplement~\cite{xiezhu2026supp}. These give $\var(\RS)\approx0.222$ and $\var(\TK)\approx0.108$.

For the fourth and sixth moments, we use Sobol quasi-Monte Carlo integration with a Brownian bridge construction \citep{caflisch1998}. The results are
$$
\E(\RS^4)\approx0.095,\qquad \E(\RS^6)\approx0.050,\qquad \E(\TK^4)\approx0.025,\qquad \E(\TK^6)\approx0.008.
$$
Together with the second moments, these values determine the polynomial densities for $m=4,5,6$ in Table~\ref{tab:polynomial-densities}. The table uses unrounded moments. Section S11 of the Supplement~\cite{xiezhu2026supp} gives the corresponding formulas.

\subsection{Finite-sample dispersion under Gaussian increments}

For each Gaussian replication, we generated two independent sequences of standard normal increments. Pearson, Spearman and Kendall correlations were computed from the paired increments and their cumulative sums. Each sample size, $n=100$, $500$ and $2000$, used 5000 independent replications, with independent runs across sample sizes.

Table~\ref{tab:moments} reports average squared correlations. Gaussian symmetry makes each coefficient's mean zero, so these estimate the corresponding variances. The level-based second moments remained close to the Pearson reference value and the two rank-limit calculations, whereas the increment-based second moments decreased with sample size. At $n=2000$, the proportions of absolute level-based correlations exceeding $0.5$ were $0.397$, $0.368$ and $0.138$ for Pearson, Spearman and Kendall, respectively. Large rank correlations therefore remained frequent despite independence of the two walks.

\begin{table}[htbp]
\caption{Empirical second moments from 5000 Gaussian replications. Parentheses give Monte Carlo standard errors. Values below $0.001$ are shown as $<0.001$.}\label{tab:moments}
\centering
\begin{tabular}{lrccc}
\toprule
Observations & $n$ & Pearson & Spearman & Kendall\\
\midrule
Random-walk levels & 100 & 0.238 (0.003) & 0.219 (0.003) & 0.108 (0.002) \\
Random-walk levels & 500 & 0.243 (0.003) & 0.224 (0.003) & 0.109 (0.002) \\
Random-walk levels & 2000 & 0.242 (0.003) & 0.224 (0.003) & 0.109 (0.002) \\
Independent increments & 100 & 0.010 ($<0.001$) & 0.010 ($<0.001$) & 0.005 ($<0.001$) \\
Independent increments & 500 & 0.002 ($<0.001$) & 0.002 ($<0.001$) & $<0.001$ ($<0.001$) \\
Independent increments & 2000 & $<0.001$ ($<0.001$) & $<0.001$ ($<0.001$) & $<0.001$ ($<0.001$) \\
\bottomrule
\end{tabular}
\end{table}

\begin{figure}[t]
\centering
\begin{minipage}[t]{0.48\textwidth}
\centering
\includegraphics[width=\linewidth]{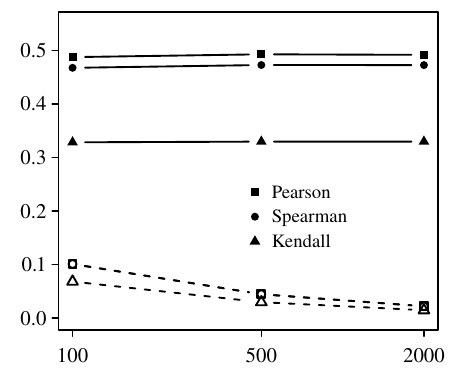}\par
{\small (a)\par}
\end{minipage}\hfill
\begin{minipage}[t]{0.48\textwidth}
\centering
\includegraphics[width=\linewidth]{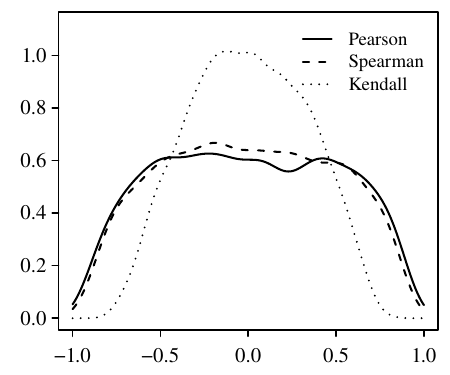}\par
{\small (b)\par}
\end{minipage}
\caption{Gaussian experiments with 5000 replications. In (a), the horizontal axis shows sample size and the vertical axis shows empirical standard deviation. Squares, circles and triangles denote Pearson, Spearman and Kendall correlations. Random-walk levels use solid lines and filled symbols, while independent increments use dashed lines and open symbols. In (b), the horizontal axis shows sample correlation and the vertical axis shows kernel density for random-walk levels at $n=2000$. Solid, dashed and dotted curves denote Pearson, Spearman and Kendall correlations.}\label{fig:dispersion}
\end{figure}

Figure~\ref{fig:dispersion}(a) shows this contrast through the standard deviations across replications. Panel (b) shows kernel density estimates of the level-based correlations at $n=2000$. These finite-sample distributions remain broad, with Kendall's coefficient less dispersed on its numerical scale.

\subsection{Calibration under Gaussian and non-Gaussian increments}

We used nominal 5\% two-sided tests. The conventional rules were the normal approximations with $Z_\rho$ from Section~\ref{sec:other-tests} for Spearman and $Z_\tau$ from Section~\ref{sec:kendall-test} for Kendall.

For Brownian calibration, the empirical 95th percentiles of the absolute coefficients in the first 2500 Gaussian replications at $n=2000$ gave thresholds of $0.813$ for Spearman and $0.597$ for Kendall. Scaled Gaussian random-walk levels have the distribution of Brownian values on the corresponding grid, so these thresholds approximate the reference quantiles in Section~\ref{sec:brownian-calibration}. We applied them unchanged to the remaining 2500 Gaussian replications and to independent non-Gaussian simulations.

The two non-Gaussian increment distributions were Student $t_3$ divided by $3^{1/2}$ and a unit-rate exponential variable minus one. Both have mean zero and variance one and satisfy Assumption~\ref{ass:random-walk}. The Student distribution has no finite third absolute moment, while the exponential distribution is asymmetric. Each setting used 2000 replications at $n=2000$. The two walks were independent in every replication.

Table~\ref{tab:calibration} reports second moments and null rejection proportions. The second moments were close to the limiting values across all three increment distributions, consistent with Theorem~\ref{thm:limit}. Conventional normal approximations rejected with proportions from $0.938$ to $0.946$, illustrating the failure described by Corollary~\ref{cor:shrinking}. The approximate Brownian thresholds gave rejection proportions from $0.047$ to $0.056$, close to the nominal level of $0.050$.

\begin{table}[htbp]
\caption{Second moments and nominal 5\% rejection proportions for random-walk levels at $n=2000$. Parentheses give Monte Carlo standard errors for second moments. The final columns use thresholds based on independent observations and the approximate Brownian reference, respectively. Evaluation used 2500 independent path pairs for Gaussian increments and 2000 for each other setting.}\label{tab:calibration}
\centering
\begin{tabular}{llccc}
\toprule
Increments & Statistic & Second moment & Independent & Brownian\\
 & & & calibration & reference\\
\midrule
Gaussian & Spearman & 0.225 (0.004) & 0.946 & 0.050 \\
Gaussian & Kendall & 0.110 (0.002) & 0.942 & 0.053 \\
Student $t_3$ & Spearman & 0.218 (0.005) & 0.941 & 0.049 \\
Student $t_3$ & Kendall & 0.105 (0.003) & 0.938 & 0.047 \\
Centred exponential & Spearman & 0.228 (0.005) & 0.946 & 0.055 \\
Centred exponential & Kendall & 0.110 (0.003) & 0.944 & 0.056 \\
\bottomrule
\end{tabular}
\end{table}

\section{Discussion}\label{sec:discussion}

Nonsense correlation persists after replacing random-walk levels by ranks. Each coefficient compares the ordering patterns of two paths at common times. Independently generated paths can have similar or opposite patterns and hence large positive or negative rank correlations. The variance formulas quantify this dispersion, while the endpoint probabilities show that both coefficients can simultaneously be arbitrarily close to either extreme.

The conventional tests fail because their critical values shrink to zero while the sample correlations remain dispersed. Marginal absolute continuity ensures that their null rejection probabilities tend to one. This also applies to the individual-observation resampling procedures in Section~\ref{sec:calibration}. Brownian quantiles control asymptotic size under Assumption~\ref{ass:random-walk}. The dependent alternative in \eqref{eq:power-alternative} shows that the calibrated tests can still have power converging to the nominal level.

The density representations characterize the limiting laws beyond their second moments and provide a counterpart to moment-based density calculations for Brownian Pearson correlation \citep{ernst2025}. Their coefficients are determined by finite Gaussian integrals, and convergence holds in $L^1$.

A further question is how to calibrate rank correlation tests for suitable classes of nonstationary time series so that they control type I error and have good power against dependent alternatives.

\section*{Supplementary material}

\noindent\emph{Supplement to ``Nonsense rank correlations between independent random walks''}. Proofs of the main results and computational details, including moment formulas, polynomial density expansions and the analysis of rank-test calibration and power.

\bibliographystyle{plainnat}
\bibliography{references,rank-additions}
\end{document}